\documentclass[12pt]{article}
\usepackage{amsmath, amscd, amssymb,latexsym, epsfig, color}
\input xy
\xyoption{all}

\newtheorem{theorem}{Theorem}[section]

\newtheorem{definition}[theorem]{Definition}

\def\proof{\smallskip\noindent {\it Proof: \ }}
\def\endproof{\hfill$\square$\medskip}
 \def\Z{\mathbb{Z}}
\def\N{\mathbb{N}}
\def\Q{\mathbb{Q}}

\def\S{\mathcal{S}}

\newcommand{\inc}{\iota}

\newcommand{\lk}{\mbox{\upshape lk}\,}
\newcommand{\cost}{\mbox{\upshape cost}\,}
\newcommand{\Hom}{\mbox{\upshape Hom}\,}
\newcommand{\field}{{\bf k}}

\newcommand{\Star}{\mbox{\upshape st}\,}

\newcommand{\soc}{\mbox{\upshape Soc}\,}
\newcommand{\Span}{\mbox{\upshape Span}\,}

\title{Gorenstein rings through face rings of manifolds.}
\author{Isabella Novik
\thanks{Research partially supported by Alfred P.~Sloan Research
Fellowship and NSF grant DMS-0500748}\\
\small Department of Mathematics, Box 354350\\[-0.8ex]
\small University of Washington, Seattle, WA 98195-4350, USA,\\[-0.8ex]
\small \texttt{novik@math.washington.edu}
\and Ed Swartz
\thanks{Research partially supported by NSF grant DMS-0600502}\\
\small Department of Mathematics, \\[-0.8ex]
\small Cornell University, Ithaca NY, 14853-4201, USA, \\[-0.8ex]
\small \texttt{ebs22@cornell.edu } }

\begin{document}
\maketitle

\begin{abstract}
The face ring of a homology manifold (without boundary) modulo a generic 
system of parameters is studied. Its socle is computed and 
it is verified that a particular quotient of this ring is Gorenstein. 
This fact is used to prove that the sphere $g$-conjecture implies all 
enumerative consequences of its far reaching generalization (due to Kalai) 
to manifolds. A special case of Kalai's manifold $g$-conjecture is 
established for homology manifolds that have a codimension-two face whose 
link contains many vertices. 
\end{abstract}

\section{Introduction}
In 1980 Stanley proved the necessity of McMullen's conjectured description 
of $f$-vectors of boundaries of simplicial convex polytopes \cite{St80}. 
At about the same time, Billera and Lee demonstrated that McMullen's 
conditions were sufficient \cite{BilleraLee}.  Since then, one of the most 
central problems in the field of face numbers of simplicial complexes is the 
$g$-conjecture. In its most optimistic form it states that, just as in the 
case of polytope boundaries,   the face ring of a homology sphere modulo a 
generic system of parameters has a Lefschetz element.  In the middle 90's 
Kalai suggested a far reaching generalization of this conjecture to all 
homology manifolds \cite[Section 7]{N98}.   It is a remarkable fact that 
all of the enumerative consequences of Kalai's conjecture are implied by 
the apparently weaker $g$-conjecture.  This follows from our main result, 
Theorem \ref{Gorenstein}, that a particular quotient of the face ring of a 
homology manifold is Gorenstein.
We also verify a special case of Kalai's conjecture when the complex 
has a codimension-two face whose link contains many vertices.

The main objects we will consider are  Buchsbaum  complexes and more 
specifically homology manifolds. Historically, Buchsbaum complexes
were defined algebraically. Here we adopt the following theorem of 
Schenzel \cite{Sch} as our definition. Let $\field$ be an infinite field 
of an arbitrary characteristic, and let $\tilde{H}_i(\Delta)$ be the
$i$-th reduced simplicial homology of $\Delta$ with coefficients in \field.

\begin{definition}
A $(d-1)$-dimensional simplicial complex $\Delta$ 
is called {\bf Buchsbaum} (over \field) if 
it is pure and for every non empty face $\tau\in \Delta$, 
$\tilde{H}_i(\lk \tau)=0$ for all $i<d-|\tau|-1$,  where
$\lk \tau =\{\sigma \in \Delta \, : \, \sigma\cap \tau=\emptyset, \,   
          \sigma\cup \tau \in \Delta\}$ is the link of $\tau$ in $\Delta$. 
\end{definition}

We say that
$\Delta$ is a {\bf homology
manifold} over $\field$ (without boundary) if it is Buchsbaum and 
in addition $\tilde{H}_{d-|\tau|-1}(\lk \tau)\cong \field$ for all 
$\emptyset\neq \tau\in \Delta$.
A \field-{\em homology sphere} is a complex $\Delta$ such that  for all 
$\tau\in \Delta,$ including $\tau=\emptyset$,
$$\tilde{H}_i(\lk \tau)=
 \left\{
\begin{array}{ll} 0 & \mbox{ if $i<d-|\tau|-1$}, \\
\field & \mbox{ if $i=d-|\tau|-1$}.
\end{array}
\right.
$$
In particular, a \field-homology sphere is a \field-homology manifold, 
and a triangulation of a topological sphere 
(topological manifold, resp.) is a
\field-homology sphere (\field-homology manifold, resp.) 
for any field \field.

If $\Delta$ is a simplicial complex on $[n]$, then
its {\em face ring} (or the {\em Stanley-Reisner ring}) is
$$\field[\Delta] := \field[x_1, \ldots, x_n]/I_\Delta,
     \quad \mbox{where }
I_\Delta=(x_{i_1}x_{i_2}\cdots x_{i_k} : \{i_1<i_2<\cdots<i_k\}
\notin\Delta).
$$
Various combinatorial and topological invariants of $\Delta$ 
are encoded in the algebraic invariants of $\field[\Delta]$ and vice versa.
For instance, if $\Delta$ is a $(d-1)$-dimensional complex, 
then  the Krull dimension of $\field[\Delta]$, $\dim\field[\Delta]$,
is equal to $d$. In this case, a set of $d$ linear forms
$\theta_1, \ldots, \theta_d \in \field[\Delta]$ is called a 
{\em linear system of parameters} (abbreviated, l.s.o.p.) if
$$\field(\Delta):=\field[\Delta]/(\theta_1, \ldots, \theta_d)$$ 
has Krull dimension zero (equivalently, $\field(\Delta)$ is a 
finite-dimensional \field-space). Assuming $\field$ is infinite,
 an l.s.o.p.~always exists:
a generic choice of $\theta_1, \ldots, \theta_d$ does the job.

One invariant that measures how far  $\Delta$ is from 
being a homology sphere is the socle of $\field(\Delta)$, 
$\soc \field(\Delta)$, where for a $\field[x_1, \ldots, x_n]$-
or $\field[\Delta]$-module $M$,
$$\soc M:= \{y\in M \, : \, x_i\cdot y =0 
\mbox{ for all } i=1,\ldots, n\}.$$ 
When $\Delta$ is a homology sphere, $\soc \field(\Delta)$
is a 1-dimensional $\field$-space.
Since $\field[\Delta]$ is a graded 
$\field$-algebra for any $\Delta$, the ring
$\field(\Delta)$ and its ideal $\soc \field(\Delta)$ are graded as well.
We denote by $\field(\Delta)_i$ and $(\soc \field(\Delta))_i$ their $i$-th
homogeneous components. It is well known (for instance, see 
\cite[Lemma III.2.4(b)]{St96}) that for any $(d-1)$-dimensional $\Delta$,
$\field(\Delta)_i$, and hence also $(\soc \field(\Delta))_i$, vanish
for all $i>d$.
If $\Delta$ is a Buchsbaum complex, then for $i\leq d$,
$(\soc(\field(\Delta))_i$ can be expressed in terms of the local cohomology
modules of $\field[\Delta]$ with respect to the irrelevant ideal, 
$H^j(\field[\Delta])$, as follows, see \cite[Theorem 2.2]{NovSw}. 

\begin{theorem} \label{socle}
Let $\Delta$ be a $(d-1)$-dimensional simplicial complex.
If $\Delta$ is Buchsbaum, then for all $0\leq i \leq d$,
  $$
(\soc \field(\Delta))_i \cong 
\left(\bigoplus_{j=0}^{d-1} \binom{d}{j}H^j (\field[\Delta])_{i-j} \right)
\bigoplus \S_{i-d},
 $$
where $\S$ is a graded submodule of $\soc H^d(\field[\Delta])$ and
$rM$ denotes the direct sum of $r$ copies of $M$.
\end{theorem}

While Theorem \ref{socle} identifies a big chunk of the socle of 
$\field(\Delta)$, its other part, $\S$, remains a mystery. 
Here we solve this mystery in the special case of a
connected {\em orientable} $\field$-homology manifold 
without boundary, thus
verifying Conjecture 7.2 of \cite{NovSw}.
A connected $\field$-homology manifold  $\Delta$ 
without boundary is called orientable if
$\tilde{H}_{d-1}(\Delta)\cong \field$. 

\begin{theorem} \label{S}
Let $\Delta$ be a $(d-1)$-dimensional connected 
 orientable \field-homology manifold without boundary. Then
$\dim_\field \S=\dim_\field \S_0=1$. In particular, 
$\dim\soc \field(\Delta)_i = \binom{d}{i}\beta_{i-1}$, where 
$\beta_{i-1}:=\dim_\field \tilde{H}_{i-1}(\Delta)$ is the $i$-th reduced
Betti number of $\Delta$.
 \end{theorem}

A graded \field-algebra of Krull dimension zero 
is called {\em Gorenstein}
if its socle is a 1-dimensional \field-space 
(see \cite[p.~50]{St96} for
many other equivalent definitions).
Let 
$$I:= \bigoplus_{i=0}^{d-1} \soc \field(\Delta)_i 
\quad \mbox{and} \quad
 \overline{\field(\Delta)}:=\field(\Delta)/I.
$$
(Note that the top-dimensional component of the 
socle is not a part of $I$). 
If $\Delta$ is a homology sphere, then $I=0$ and 
$\overline{\field(\Delta)}=\field(\Delta)$ is a Gorenstein ring 
\cite[Theorem II.5.1]{St96}.
What if $\Delta$ is a homology manifold other than a sphere? How far is
$\field(\Delta)$ from being Gorenstein in this case? 
The answer (that was conjectured in \cite[Conjecture 7.3]{NovSw})
turns out to be surprisingly simple:

\begin{theorem} \label{Gorenstein}
Let $\Delta$ be a $(d-1)$-dimensional connected simplicial complex.
If $\Delta$ is an orientable \field-homology manifold without boundary, then
$\overline{\field(\Delta)}$ is Gorenstein.
\end{theorem}

In \cite{Sch}, Schenzel computed the Hilbert function of $\field(\Delta)$
for a Buchsbaum complex $\Delta$ in terms of its face and
Betti numbers. It follows from Theorem \ref{S}
combined with Schenzel's formula and Dehn-Sommerville relations
\cite{Klee} that for a connected
orientable homology manifold $\Delta$, the Hilbert function of
$\overline{\field(\Delta)}$ is symmetric, that is, 
\begin{equation}  \label{symmetry}
\dim_\field \overline{\field(\Delta)}_i = 
\dim_\field \overline{\field(\Delta)}_{d-i} \quad \mbox{ for all } 
0\leq i \leq d.
\end{equation}
As the Hilbert function of a Gorenstein ring of Krull dimension zero
is always symmetric \cite[p.~50]{St96},
Theorem \ref{Gorenstein} gives an alternative
algebraic proof of (\ref{symmetry}).

Theorems \ref{S} and \ref{Gorenstein} are ultimately related 
to the celebrated $g$-conjecture
that provides a complete characterization of possible face numbers 
of homology spheres. The most optimistic version of this conjecture
is a very strong manifistation of the symmetry of the Hilbert function. 
It asserts that if $\Delta$ is a $(d-1)$-dimensional homology sphere, 
and $\omega$ and $\Theta=\{\theta_1, \ldots, \theta_d\}$ 
are sufficiently generic linear forms, 
then for $i\leq d/2$, multiplication
$$
\omega^{d-2i}: \field(\Delta)_i \longrightarrow \field(\Delta)_{d-i}
$$
is an isomorphism. At present this conjecture is known to hold
only for the class of polytopal spheres \cite{St80} and edge 
decomposable spheres \cite{Murai}, \cite{Nevo}.

Kalai's far-reaching generalization 
of the $g$-conjecture \cite{N98} 
posits that if $\Delta$ is an orientable
homology manifold  and 
$\omega, \Theta$  are sufficiently generic, then
$$
\omega^{d-2i}: \overline{\field(\Delta)}_i \longrightarrow 
\overline{\field(\Delta)}_{d-i}
$$ 
is still an isomorphism for all $i\leq d/2$.  
Let $h''_i = \dim_\field \overline{\field(\Delta)}_i.$ 
Given a system of parameters and $\omega$ which satisfy Kalai's conjecture 
it is immediate that multiplication 
$\omega:  \overline{\field(\Delta)}_i \to 
    \overline{\field(\Delta)}_{i+1}$ is an injection 
for $i < d/2.$  So, $h''_0 \le \dots \le h''_{\lfloor d/2 \rfloor},$ 
and examination of $ \overline{\field(\Delta)} / (\omega)$ 
shows that the nonnegative integer vector 
$(h''_0, h''_1 - h''_0, \dots, 
h''_{\lfloor d/2 \rfloor} - h''_{\lfloor d/2 \rfloor - 1})$ is an M-vector, 
i.e.  satisfies Macaulay's nonlinear arithmetic conditions 
(see \cite[p.~56]{St96}) for the Hilbert 
series of a homogeneous quotient of a polynomial ring.  
Applying the same reasoning to $I_{[i]}:=\oplus_{j=0}^{i} I_j$, which is
also an ideal, and $\field(\Delta, i):= \field(\Delta)/I_{[i]}$, instead of
$I$ and $\overline{\field(\Delta)} = \field(\Delta)/I$,  the
second conclusion can be strengthen to $(h''_0, h''_1-h''_0,
\ldots, h''_i-h''_{i-1}, h''_{i+1}-h''_{i}+\binom{d}{i+1}\beta_i)$ is an
$M$-vector for every $i<\lfloor d/2\rfloor$.
In fact, as we will see in Theorem~\ref{connection}, 
these two conclusions follow from the $g$-conjecture.

For Kalai's conjecture we prove this special case. 

\begin{theorem} \label{g-thm}
Let $\Delta$ be a $(d-1)$-dimensional
orientable \field-homology manifold with $d\geq 3$. 
If $\Delta$ has a $(d-3)$-dimensional face $\tau$ whose link contains
all of the 
vertices of $\Delta$ that are not in $\tau$, then for generic choices of
$\omega$ and $\Theta$,
$\omega^{d-2}: \overline{\field(\Delta)}_1 \longrightarrow 
\overline{\field(\Delta)}_{d-1}
$
is an isomorphism.
\end{theorem}
The condition that the link of $\tau$ contains all of the vertices of
$\Delta$ that are not in $\tau$ is equivalent to saying that  
every vertex of $\Delta$ is in the star of $\tau$,
$\Star \tau:=\{\sigma\in \Delta \, : \, \sigma\cup \tau\in \Delta\}$. 
This condition is not as restrictive as one might think:
the results of \cite[Section 5]{Sw} imply that every connected 
homology manifold $M$ without boundary that has a triangulation, 
always has a triangulation $\Delta$ satisfying this condition.

The outline of the paper is as follows. 
We verify Theorem \ref{S} in Section 2. 
The main ingredient in the proof is  Gr\"abe's 
explicit description of  $H^d(\field[\Delta])$
as a $\field[\Delta]$-module in terms of the simplicial
(co)homology of the links of faces of $\Delta$ and maps between
them \cite{Grabe}. Theorem~\ref{Gorenstein} is proved in Section 3
and is used to explore the relationship between 
the $g$-conjecture and Kalai's conjecture.
Section 4 is devoted to the proof of Theorem~\ref{g-thm}. The proofs of
Theorems \ref{Gorenstein} and \ref{g-thm}
rely heavily on a result from \cite{Sw} relating $\field(\lk v)$ 
(for a vertex $v$ of $\Delta$) to the
principal ideal $(x_v) \subset \field(\Delta)$.

\section{Socles of homology manifolds}
The goal of this section is to verify Theorem \ref{S}. 
To do so we analyze $\soc H^d(\field[\Delta])$ and prove the following.

\begin{theorem} \label{socH^d}
If $\Delta$ is a connected orientable $(d-1)$-dimensional
\field-homology manifold without boundary, then
$\soc H^d(\field[\Delta])_i=0$ for all $i\neq 0$.
\end{theorem}

The proof relies on results from \cite{Grabe}.
We denote by $|\Delta|$ the geometric realization
of $\Delta$. For a face $\tau\in \Delta$, 
let $\cost \tau :=\{\sigma \in \Delta \, : \sigma\not\supset \tau\}$
be the contrastar of $\tau$,
let $H^i(\Delta, \cost \tau)$ be the simplicial $i$-th
cohomology of a pair (with coefficients in \field), and for
 $\tau \subset \sigma \in \Delta$, let $\inc^*$ be the map
$H^i(\Delta, \cost \sigma) \to H^i(\Delta, \cost \tau)$ induced
by inclusion $\inc: \cost\tau \to \cost\sigma$. 
Also, if $\emptyset \neq \tau\in \Delta$, 
let $\hat{\tau}$ be the barycenter of $\tau$.
Finally, for a vector $U=(u_1, \ldots, u_n)\in \Z^n$, let 
$s(U):=\{l : u_l\neq 0\}\subseteq [n]$ be the support of $U$, 
let $\{e_l\}_{l=1}^n$ be the standard basis for $\Z^n$, and let $\N$ 
denote the set of nonnegative integers. 

\begin{theorem} {\rm {\bf [Gr\"abe]}}   \label{Grabe}
The following is an isomorphism of $\Z^n$-graded $\field[\Delta]$-modules 
\begin{equation}  \label{Hochster}
H^{i+1}(\field[\Delta]) \cong 
\bigoplus_{\mbox{\tiny${\begin{array}{cc} -U\in \N^n\\
                             s(U)\in\Delta
           \end{array}}$}}
H^i (\Delta, \cost s(U)),
\end{equation}
where the $\field[\Delta]$-structure on the $U$-th component
of the right-hand side is given by 
$$
\cdot x_l = \left\{ \begin{array}{lll}
\text{$0$-map}, & \mbox{ if } l\notin s(U)\\
\text{identity map}, & \mbox{ if } l\in s(U) \mbox{ and } l\in s(U+e_l)\\
\inc^*: H^{i}(\Delta, \cost s(U)) \to H^i(\Delta, \cost s(U+e_l)), &
\mbox{ otherwise}. 
\end{array}
\right.
$$
\end{theorem}
We remark that the isomorphism of (\ref{Hochster}) 
on the level of vector spaces
(rather than $\field[\Delta]$-modules) is due to Hochster, see 
\cite[Section II.4]{St96}.

\medskip\noindent {\it Proof of Theorem \ref{socH^d}: \ }
In view of Theorem \ref{Grabe}, to prove that 
$\soc H^d(\field[\Delta])_i=0$
for all $i\neq 0$, it is enough to show that for every 
$\emptyset \neq \tau\in\Delta$ and $l\in\tau$, the map
$\inc^*: H^{d-1}(\Delta, \cost \tau) \to 
H^{d-1}(\Delta, \cost \sigma)$, where
$\sigma= \tau-\{l\}$, is an isomorphism. 
Assume first that $\sigma \neq \emptyset$.
Consider the following diagram.
$$
\begin{CD}
H^{d-1}(|\Delta|) @>(j^\star)^{-1}>> 
H^{d-1}(|\Delta|, |\Delta| - \hat{\sigma})
 @>f^\star>>  H^{d-1}(\Delta, \cost \sigma) \\
 @|  @. @A\inc^\star AA \\
 H^{d-1}(|\Delta|) @>(j^\star)^{-1}>> 
H^{d-1}(|\Delta|, |\Delta| - \hat{\tau}) 
@>f^\star>>  H^{d-1}(\Delta, \cost \tau) 
\end{CD}
$$
The two $f^\star$ maps are induced by inclusion and are 
isomorphisms by the usual deformation retractions. 
$j^\star$ is also induced
by inclusion. Since $\Delta$ is connected and
orientable, all of the spaces are one-dimensional and $j^\star$ 
is an isomorphism, so that $(j^\star)^{-1}$ is well-defined and is 
an isomomorphism as well. Hence compositions
$f^\star \circ (j^\star)^{-1}$ are isomorphisms.
The naturality of $j^\star$ implies
that the diagram is commutative.   
It follows that $\inc^\star$ is an isomorphism. 
If $\sigma=\emptyset$, replace in the above diagram
$H^{d-1}(|\Delta|, |\Delta| - \hat{\sigma})$ with 
$H^{d-1}(|\Delta|)$. The same reasoning applies.
\endproof
 
We are now in a position to complete the proof of Theorem \ref{S}.

\medskip\noindent {\it Proof of Theorem \ref{S}: \ }
Since 
$\S_i$ is a subspace of $\soc H^d(\field[\Delta])_i$
(Theorem~\ref{socle}), and since the later space is the zero-space
whenever $i\neq 0$ (see Theorem  \ref{socH^d}), it follows that $\S_i=0$
for all $i\neq 0$. For $i=0$, we have 
$$\dim \S_0=\dim \soc \field(\Delta)_d= \dim \field(\Delta)_d=
\beta_{d-1}(\Delta)=1.$$
Here the first step is by Theorem \ref{socle}, 
the second step is a consequence
of $\field(\Delta)_d$ being the last 
non vanishing homogeneous component
of $\field(\Delta)$, and the third 
step is by Schenzel's formula \cite{Sch}
(see also \cite[Theorem II.8.2]{St96}). 
The ``In particular''-part then follows
from Theorem~\ref{socle}, isomorphism (\ref{Hochster}), 
and the standard fact that 
$H^{i-1}(\Delta, \cost \tau) \cong \tilde{H}_{i-|\tau|-1}(\lk \tau)$,
see e.g.~\cite[Lemma 1.3]{Grabe}.
\endproof

\section{Gorenstein property}
In this section we verify Theorem \ref{Gorenstein} and use it to discuss
connections between various $g$-conjectures. To prove that 
$\overline{\field(\Delta)}=\field(\Delta)/I$ is Gorenstein, 
where $\Delta$ is a  $(d-1)$-dimensional connected
orientable homology manifold and $I=\oplus_{i=0}^{d-1} \soc \field(\Delta)_i$,
we have to check that the operation of 
moding out by $I$ does not introduce new socle elements.
This turns out to be a simple application of \cite[Proposition 4.24]{Sw},
which we review now.

Let $\Theta=\{\theta_1, \ldots, \theta_d\}$ 
be an l.s.o.p.~for $\field[\Delta]$
and let $v$ be a vertex of $\Delta$.
Fix a facet $\tau=\{v=v_1, v_2, \ldots, v_d\}$ that contains $v$. 
By doing Gaussian elimination on the $d\times n$-matrix whose $(i,j)$-th
entry is the coefficient of $x_j$ in $\theta_i$ we can
assume without loss of generality that 
$\theta_i=x_{v_i}+\sum_{j\notin\tau}\theta_{i,j}x_j$. 
Denote by $\theta'_i$ the linear form obtained from $\theta_i$ by removing 
all summands involving $x_j$ for $\{j\}\notin \lk v$. 
Then
$\Theta':=\{\theta'_2, \ldots, \theta'_d\}$ can be considered as 
a subset of $\field[\lk v]_1$. Moreover,
 it is easy to check, say, using \cite[Lemma III.2.4(a)]{St96},
that $\Theta'$ forms an l.s.o.p.~for $\field[\lk v]$. The ring
$\field(\lk v):=\field[\lk v]/(\Theta')$ has a natural 
$\field[x_1,\ldots,\hat{x_v}, \ldots, x_n]$-module structure (if
$j\neq v$ is not in the link of $v$, then multiplication by $x_j$
is the zero map), and defining
$$x_v \cdot y:=-\theta'_1 \cdot y \quad \mbox{ for } y\in \field(\lk v)$$ 
extends it to a $\field[x_1, \ldots, x_n]$-module structure. 
Proposition 4.24 of \cite{Sw} asserts the following.

\begin{theorem}   \label{link}
Let $\Delta$ be an orientable homology manifold.
The map 
$$\phi: \field[\lk v]/(\Theta') 
\to (x_v)\left(\field[\Delta]/(\Theta)\right)
\quad 
\mbox{given by }z \mapsto x_v\cdot z,$$
is well-defined and is 
an isomorphism (of degree 1) of $\field[x_1, \ldots, x_n]$-modules.
Its inverse, $x_v \cdot z \mapsto\overline{z}$, is given by replacing each 
occurrence of $x_v$ in $z$ with $-\theta'_1$ and setting all $x_j$ for
$j\neq v$ not in the link of $v$ to zero.
\end{theorem}

\medskip\noindent {\it Proof of Theorem \ref{Gorenstein}: \ }
To prove the theorem it is enough to show that the socle of
$\overline{\field(\Delta)}=\field(\Delta)/I$, where 
$I=\bigoplus_{j=0}^{d-1} \soc \field(\Delta)_j$, vanishes in all degrees
$j\neq d$. This is clear for $j=d-1$. For $j\leq d-2$,
consider any element $y\in \field(\Delta)_j$ 
such that $x_v\cdot y\in \soc \field(\Delta)$ for all $v\in [n]$. 
We have to check that $y\in \soc \field(\Delta)$. 
And indeed, the isomorphism of Theorem \ref{link} implies that
$\overline{y}:=\phi^{-1}(x_v\cdot y)\in \field(\lk v)_j$ is in the socle of
$\field(\lk v)$. Since $\lk v$ is a $(d-2)$-dimensional
homology sphere, $\field(\lk v)$ is Gorenstein, 
and hence its socle vanishes in all degrees
except $(d-1)$-st one. Therefore, $\overline{y}=0$, and hence 
$x_v\cdot y=\phi(\overline{y})=0$  in $\field(\Delta)$.
Since this happens for all $v\in [n]$, it follows that 
$y\in \soc \field(\Delta)$. \endproof

We now turn to discussing the sphere and manifold $g$-conjectures 
and the connection between them. As was mentioned in the introduction, 
the strongest $g$-conjecture for homology spheres and its 
generalization (due to Kalai) for homology manifolds asserts that
if $\Delta$ is a $(d-1)$-dimensional connected orientable 
homology manifold, then for generically chosen $\omega$ and 
$\Theta=\{\theta_1, \ldots, \theta_d\}$ in $\field[\Delta]_1$,
the map
 $$\cdot \omega^{d-2i}: \overline{\field(\Delta)}_i  \to 
\overline{\field(\Delta)}_{d-i}
\mbox{ is an isomorphism for all $i\leq \lfloor d/2\rfloor$}.
$$ 
We refer to this conjecture as the {\em  strong (sphere or manifold)
$g$-conjecture}.
If true, it would imply that 
$\cdot \omega: \overline{\field(\Delta)}_i  \to 
\overline{\field(\Delta)}_{i+1}
$
is injective for all $i < \lfloor d/2\rfloor$ and is surjective for
all $i \geq \lceil d/2\rceil$. We refer to this weaker statement as the
{\em (sphere or manifold) $g$-conjecture}. 
(Both, the stronger and the weaker conjectures yield exactly the same
combinatorial restrictions on the face numbers of $\Delta$.)
Clearly, the manifold $g$-conjecture 
implies the sphere $g$-conjecture. 
The following result shows that they are almost equivalent: 
the strong sphere $g$-conjecture in the middle degree
implies the manifold $g$-conjecture in all degrees.

\begin{theorem}  \label{connection}
Let $\Delta$ be a $(d-1)$-dimensional connected 
orientable homology manifold. If for at least $(n-d)$ of the  
vertices $v$ of $\Delta$
and generically chosen $\omega$ and $\Theta'$ in $\field[\lk v]_1$, the map
$\cdot \omega: \field(\lk v)_{\lfloor (d-1)/2\rfloor}  \to 
\field(\lk v)_{\lfloor(d-1)/2\rfloor+1}$ is surjective, then
$\Delta$ satisfies the manifold $g$-conjecture.
\end{theorem}

\proof The condition on the links implies by \cite[Theorem 4.26]{Sw} 
that for generic choices of $\omega$ and $\Theta$ in $\field[\Delta]_1$,
the map $\cdot \omega: \field(\Delta)_{\lceil d/2\rceil}
\to \field(\Delta)_{\lceil d/2\rceil+1}$ is surjective. 
Hence the map 
$\cdot \omega: \overline{\field(\Delta)}_{\lceil d/2\rceil}
\to \overline{\field(\Delta)}_{\lceil d/2\rceil+1}$ is 
surjective. Thus, $\overline{\field(\Delta)}/(\omega)_i$ 
vanishes for $i=\lceil d/2\rceil+1$, and hence also for all 
$i>\lceil d/2\rceil+1$, and we infer that 
$\cdot \omega: \overline{\field(\Delta)}_{i}
\to \overline{\field(\Delta)}_{i+1}$ is surjective
for $i\geq \lceil d/2\rceil$. This in turn yields that the dual map
$\cdot\omega: \Hom_\field(\overline{\field(\Delta)}_{i+1}, \field) \to 
\Hom_\field(\overline{\field(\Delta)}_{i}, \field)$ is injective for
all $i\geq \lceil d/2\rceil$. Since  $\overline{\field(\Delta)}$ is 
Gorenstein,
$\Hom_\field(\overline{\field(\Delta)}_{i}, \field)$ is 
naturally isomorphic to $\overline{\field(\Delta)}_{d-i}$
(see Theorems I.12.5 and I.12.10 in \cite{St96}).
Therefore, $\cdot \omega: \overline{\field(\Delta)}_{j}
\to \overline{\field(\Delta)}_{j+1}$ is injective for 
$j<\lfloor d/2 \rfloor$. \endproof

Theorem \ref{connection} combined with Stanley's $g$-theorem for polytopes
\cite{St80}, implies that every $(d-1)$-dimensional connected orientable
$\Q$-homology manifold all of whose vertex links
are polytopal spheres satisfies the manifold $g$-conjecture.

\section{A special case of the $g$-theorem}
In this section we prove Theorem \ref{g-thm}. 
As in the proof of Theorem~\ref{Gorenstein} we will
rely on Theorem~\ref{link} and notation introduced there. 
Since the set of all $(\omega, \Theta)$
for which $\cdot \omega^{d-2}: \overline{\field(\Delta)}_1 \to
\overline{\field(\Delta)}_{d-1}$ is an isomorphism, is
 a Zariski open set (see \cite[Section 4]{Sw}), 
it is enough to find one $\omega$ that does the job. Surprisingly, the 
$\omega$ we find is ``very non-generic": as we will see, 
$w=x_v$ where $v\in \tau$ does the job.

\medskip\noindent {\it Proof of Theorem \ref{g-thm}: \ }
We use induction on $d$ starting with $d=3$. 
In this case, the face $\tau$ is simply a vertex, say $v$.
Let $\sigma=\{v=v_1, v_2, v_3\}$ be a facet containing $v$.
Let $\Theta$ be a generic l.s.o.p.~for $\field[\Delta]$, and
as in Section~3 assume that 
$\theta_i=x_{v_i}+\sum_{j\notin\sigma}\theta_{i,j}x_j$.
Since $\dim_\field \overline{\field(\Delta)}_1=
\dim_\field \overline{\field(\Delta)}_2$ (see Eq.~(\ref{symmetry}))
and since $\dim\soc \field(\Delta)_1=0$ (e.g.~by Theorem~\ref{S}), 
we will be done if we check that the map
$\cdot x_{v}:  \field(\Delta)_1 \to \overline{\field(\Delta)}_2$
is injective. And indeed, if $x_{v}\cdot z \in \soc \field(\Delta)_2$
for some $z\in \field(\Delta)_1$, then by 
the isomorphism of Theorem~\ref{link}, 
$\overline{z}=\phi^{-1}(x_{v}\cdot z)$ is in $\soc \field(\lk v)_1$.
But the later space is the zero space, so $\overline{z}=0$ in 
$\field(\lk v)_1$, that is, $\overline{z}\in \Span (\theta'_2, \theta'_3)$. 
Since $\lk v$ contains all the vertices of $\Delta$ except $v$, 
it follows that  $\theta'_2=\theta_2$, $\theta'_3=\theta_3$,
and $z-\overline{z}$ is a multiple of 
$\theta_1$, and hence that $z=0$ in $\field(\Delta)$.

The proof of the induction step goes along the same lines: 
let $\tau=\{v_1, \ldots, v_{d-2}\}$, 
let $\sigma=\tau\cup \{v_{d-1}, v_d\}$ be any facet containing $\tau$,
and let $\Theta$ be a generic ls.o.p.~for $\field[\Delta]$.
We show that for $v\in\tau$, say $v=v_1$, the map 
$\cdot x_{v}^{d-2}:  \field(\Delta)_1 \to 
\overline{\field(\Delta)}_{d-1}$ is an injection. If 
$x_{v}^{d-2}\cdot z \in \soc \field(\Delta)_{d-1}$
for some $z\in \field(\Delta)_1$, then by Theorem~\ref{link},
$$x_{v}^{d-3}\cdot \overline{z}=
(-\theta'_1)^{d-3}\cdot \overline{z} \in \soc \field(\lk v)_{d-2}.$$
However, $\lk v$ is a $(d-2)$-dimensional homology sphere satisfying
the same assumptions as $\Delta$: the star of the face $\tau-\{v\}$ 
contains all the vertices of $\lk v$. Hence by inductive hypothesis
applied to $\lk v$ with $\omega=-\theta'_1$ and 
$\Theta'=(\theta'_2,\ldots, \theta'_d)$, 
the map $\cdot (-\theta'_1)^{d-3}: 
\field(\lk v)_1 \to \overline{\field(\lk v)}_{d-2}$ is an injection.
Thus $\overline{z}=0$ in $\field(\lk v)_1$, 
and so $z=0$ in $\field(\Delta)_1$.
\endproof

\smallskip\noindent{\bf Remark} \;\;
Barnette's lower bound theorem \cite{Ka87} asserts that
for every $d$ and $n$, the $f$-vector
of any $(d-1)$-dimensional connected triangulated
manifold with $n$ vertices is
minimized componentwise by the $f$-vector of
(the boundary of) a stacked $d$-polytope with
$n$ vertices, $S(d,n)$. There is a conjectural algebraic
strengthening of this result  based on Kalai's notion of
{\em algebraic shifting}, an operation that was
introduced in the mid-eighties as a tool for studying $f$-numbers of
simplicial complexes (see e.g.~\cite{BjKa} and a survey paper
\cite{Ka02}). This conjecture posits that the algebraic shifting
of every connected orientable $(d-1)$-dimensional manifold on
$n$ vertices contains as a subset $\Delta(S(d,n))$ ---
the algebraic shifting of $S(d,n)$.
The complex  $\Delta(S(d,n))$  was recently computed by
Nevo \cite[Example 2.18]{Nevo} and independently by Murai \cite{Mur07}.
Using their result, together with standard facts on
rev-lex generic initial ideals \cite[Prop.~15.12]{Eis}
and an algebraic definition
of Buchsbaum complexes \cite[Def.~and Thm.~3.1(ii)]{Sch},
it is easy to show that in the case of
symmetric algebraic shifting,
a $(d-1)$-dimensional connected orientable
$\Q$-homology manifold $\Delta$ satisfies the above conjecture
if and only if for generically
chosen $\theta_1,\ldots, \theta_d$ and $\omega$, the map
$\cdot \omega^{d-2}: \overline{\field(\Delta)}_1 \to
\overline{\field(\Delta)}_{d-1}$ is a $\Q$-isomorphism.
In particular, it follows that this conjecture holds for
every $(d-1)$-dimensional orientable $\Q$-homology
manifold satisfying assumptions of Theorem \ref{g-thm}.

\end{document}